\documentclass[preprint,10pt]{elsarticle}

\usepackage{amssymb,amsmath}
\def\be{\begin{equation}}
\def\ee{\end{equation}}

\usepackage{graphicx,eucal}
\usepackage[all]{xy}
\xyoption{all}

\usepackage{color}

\journal{Physica A}

\begin{document}

\begin{frontmatter}



\title{Immortal Branching Processes}
\author{P. L. Krapivsky}
\ead{paulk@bu.edu}
\address{Department of Physics, 590 Commonwealth
  Avenue, Boston University, Boston, MA 02215, USA}
\author{S. Redner}
\ead{redner@santafe.edu}
\address{Santa Fe Institute, 1399 Hyde Park Road, Santa Fe, New Mexico 87501,
  USA}

\begin{abstract}
  We investigate the dynamics of an \emph{immortal} branching process.  In
  the classic, critical branching process, particles give birth to a single
  offspring or die at the same rates.  Even though the average population is
  constant in time, the ultimate fate of the population is extinction.  We
  augment this branching process with immortality by positing that either:
  (a) a single particle cannot die, or (b) there exists an immortal stem cell
  that gives birth to ordinary cells that can subsequently undergo critical
  branching.  We also treat immortal two-type branching.  We discuss the new
  dynamical aspects of these immortal branching processes.
\end{abstract}

\end{frontmatter}

\section{Introduction}

We dedicate this contribution to the memory of Dietrich Stauffer, a singular
personality in the physics community.  He was a man with strong views, who
spoke with few compunctions, and who had a great sense of humor.  Looking
back at the emails that he exchanged with one of us (SR) brought back
pleasant memories of his caustic wit.  One such example is Dietrich's email
to another person (with a cc to me) in which he writes: ``I am very grateful
that you try to educate the Canadian citizen Sid Redner [yes, SR is Canadian]
in a civilised language. But since he is a traitor working for the Southern
Enemy [meaning the USA], I doubt that you will have a long-lasting
success. He does not even use the spelling of Her Britannic Majesty, Canada's
head of state.''  He enjoyed making fun of SR's allegedly traitorous behavior
(or behaviour) and SR enjoyed getting poked by him.  In addition to the many
verbal and written barbs that we exchanged, it was also a great pleasure to
discuss science and indeed any field of academic study with him.  He will be
sorely missed.

Dietrich is best known for his contributions to the field of
percolation~\cite{stauffer}.  Percolation has applications to diverse fields;
one such example is epidemic spread.  From the dynamical perspective, the
temporal course of an epidemic can be viewed as a birth/death or a branching
process.  Namely, an infected individual can infect someone who is not yet
sick and an infected individual can be ``removed'' by either recovering or
expiring.  If we denote an infected individual as P, the epidemic dynamics
can be symbolically represented as
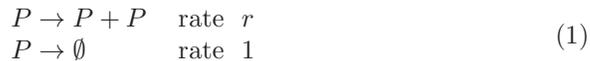
\begin{equation}
\label{modeldef}
 \begin{tabular}{ll}
 $P \to  P+P$                & rate $~r$\\
 $P \to  \emptyset$       & rate $~1$\\
\end{tabular}
\end{equation}
where the death rate is set to one by choosing appropriate time units.  The
\emph{critical} branching process in which the birth and death rates are the
same, $r=1$, is especially popular, and has numerous
applications~\cite{feller,teh,vatutin}.  Although the adjective ``critical''
suggests that this branching process is idiosyncratic because of the tuning
of the birth and death rates, this case is arguably more important than its
subcritical and supercritical brethren.  Indeed, subcritical branching
($r<1$) quickly goes extinct while supercritical branching ($r>1$) results in
unlimited growth.

Branching processes are related to the stochastic version of the
susceptible-infected-recovered (SIR) infection model, in which Ps correspond
to infected individuals and $\emptyset$s to recovered.
In the SIR model~\cite{ntjb,ma,antia,lux}, the population consists of
susceptible, infected, and recovered individuals.  The infection spreads by
contact between infected and susceptible individuals, while infected
individuals spontaneously recover or die.  In the context of infectious
diseases, critical branching processes arise as a balance between human
efforts that strive to reduce the infection rate~\cite{ma} and natural
evolution that increases the infectivity of diseases~\cite{antia}.  Critical
branching can also describe homeostasis (the maintenance of healthy
tissues)~\cite{clayton07,JS08}.  For instance, the
maintenance of skin cells arises when progenitor cells (P), which exist far
below the epidermis, undergo critical branching and (instead of death)
differentiate into post-mitotic cells that are eventually shed.

In this contribution, we study \emph{immortal} branching processes (IBPs).
Our first example is identical to the critical branching process with a
simple twist---when a single cell remains, it cannot die.  One of the
motivations for this type of modeling is the long-term dynamics of HIV
infection in which a nearly undetectable level of viral load can linger in
patients for a decade or longer and then subsequently increase~\cite{hiv}.
This branching without the possibility of extinction greatly oversimplifies
the HIV dynamics, but perhaps it points the way to more realistic models.  We
analyze this model in Sec.~\ref{sec:no} and demonstrate unexpectedly subtle
dynamics that are replete with logarithmic corrections in the long time
limit.

Another IBP is a critical branching process that is augmented by the presence
of a single immortal stem cell whose properties are not affected by
progenitor cells.  The reaction steps of this IBP are
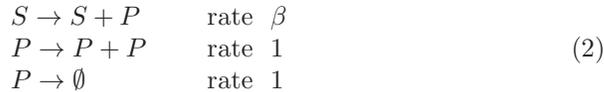
\begin{equation}
\label{modeldefi}
 \begin{tabular}{ll}
 $S \to  S+P$ ~~~         & rate $~\beta$\\
 $P \to  P+P$                & rate $~1$\\
 $P \to  \emptyset$       & rate $~1$\\
 \end{tabular}
\end{equation}
The first process accounts for a stem cell (S) that creates a single
progenitor cell P, while the latter two reactions represent critical
branching.  Mathematically, the reactions \eqref{modeldefi} correspond to
critical branching with input.  The strength of the source, which is
quantified by $\beta$, can vary with time in principle.  For simplicity, we
focus on the case where $\beta=\text{const}$.  There is no extinction because
the stem cell is immortal and the average number of progenitor cells grows as
$\beta t$.

\section{Critical Branching}
\label{sec:classic}

It is helpful to first recall basic facts about critical branching.  The
probability $P_m(t)$ to have $m$ cells at time $t$ satisfies the rate
equation
\begin{equation}
\label{master}
\frac{dP_m}{dt} = (m+1)P_{m+1} -2mP_m+(m-1)P_{m-1}\,.
\end{equation}
The first term on the right-hand side corresponds to any one of $m+1$ cells
dying with rate 1; the other terms have similar explanations.  We always
assume that the system initially contains a single cell,
$P_m(t=0) = \delta_{m,1}$.

There are various ways to solve Eqs.~\eqref{master}. Perhaps the simplest
makes use of the exponential ansatz (see, e.g.,~\cite{book})
\begin{equation}
\label{exp}
P_m(t) =A(t) [a(t)]^{m-1},
\end{equation}
which often works for infinite sets of equations with linear dependences on
$m$.  With ansatz \eqref{exp}, the infinite set of equations \eqref{master}
reduces to $\frac{dA}{dt}=2A(a-1)$ and $\frac{da}{dt} =(a-1)^2$, from which
$A=(1+t)^{-2}$ and $a=t/(1+t)$.  Thus
\begin{equation}
\label{classical}
P_m(t) = \frac{t^{m-1}}{(1+t)^{m+1}}
\end{equation}
is the probability that there are $m$ cells at time $t$. To appreciate the
meaning of this solution, it is useful to look at the integer moments
\begin{equation*}
%\label{m-av}
\langle m^k \rangle = \sum_{m\geq 1}m^k P_m(t)\,.
\end{equation*}
The first few moments are
$\langle m^0 \rangle = (1+t)^{-1}, ~\langle m^1 \rangle = 1, ~\langle m^2
\rangle = 1+2t$.  The zeroth moment is just the probability that the
population does not go extinct---the survival probability.  Even though
$\langle m^0\rangle \to 0$ as $t\to\infty$ so that eventual extinction
necessarily occurs, the mean extinction time is infinite.  The conserved
nature of the first moment reflects the equality of the birth and death rates
so that the average population size is constant.  In contrast, for those
realizations that do not go extinct, the population grows linearly with time,
as manifested by the behavior of the second moment.

\section{Branching Without Extinction}
\label{sec:no}

We now amend critical branching so that extinction is impossible.  A simple
way to impose this constraint is to retain the rules of critical branching
but forbid death when only a single cell remains.  That is,
Eqs.~\eqref{master} remain valid for $m\geq 2$, while the probability
$P_1(t)$ now evolves according to
\begin{equation}
\label{single}
\frac{dP_1}{dt} =  -P_1+2P_2\,.
\end{equation}

To solve Eqs.~\eqref{master} and \eqref{single}, first notice
that the ansatz \eqref{exp} no longer solves them, which suggests that the
solution does not have an exponential form.  Instead, we proceed
conventionally and introduce the Laplace transform
\begin{equation*}
\widetilde{P}_m(s) = \int_0^\infty dt\,e^{-st} P_m(t)\,,
\end{equation*} 
to recast \eqref{master} and \eqref{single} into the linear algebraic
equations
\begin{align}
\label{Lap:2}
\begin{split}
  (s+2m)\widetilde{P}_m &=
(m-1)\widetilde{P}_{m-1}+(m+1)\widetilde{P}_{m+1}\qquad m\geq 2\,,\\
(s+1)\widetilde{P}_1 &= 2\widetilde{P}_2 + 1\,.
\end{split}
\end{align}

We now apply the generating function method.  By standard steps, we find that
the generating function
\begin{equation*}
%\label{gen}
\mathcal{P}(s,z) = \sum_{m\geq 1} \widetilde{P}_m(s) \, z^m
\end{equation*}
satisfies the ordinary differential equation (ODE)
\begin{equation}
\label{gen:eq}
s\mathcal{P} = (z-1)^2\, \frac{d \mathcal{P}}{dz} + (z-1)\widetilde{P}_1(s) - z\,.
\end{equation}
The details of the solution to this equation and the inversion of the Laplace
transform are given in \ref{app:a}.  The final result for the population
distribution in the Laplace domain is
\begin{equation}
\label{Pm:final} 
\widetilde{P}_m(s) = \frac{1}{s e^s \Gamma(0,s)}\int_0^\infty d\eta\,e^{-s\eta}\,\frac{\eta^{m-1}}{(1+\eta)^{m+1}}\,.
\end{equation}

While it does not seem possible to invert this Laplace transform and obtain
$P_m(t)$ in terms of elementary or standard special functions, we can
determine the long-time behavior from the $s\to 0$ asymptotic of the Laplace
transform.  Using the asymptotic formula (see, e.g.,~\cite{BO78})
\begin{equation*}
\Gamma(0,s) = -\ln s - \gamma + \sum_{k\geq 1}\frac{(-s)^k}{k\times k!}\,,
\end{equation*}
where $\gamma=0.577\,215\ldots$ is Euler's constant and
$\Gamma(0,s) = \int_s^\infty du\,e^{-u}/u$ is the incomplete Gamma function,
the leading $s\to 0$ behavior of $\widetilde{P}_1(s)$ is
\begin{subequations}
\begin{equation}
\label{P1:asymp-Lap}
\widetilde{P}_1(s) \simeq \frac{1}{s(-\ln s - \gamma)} \,,
\end{equation}
from which we deduce the unexpectedly slow long-time decay
\begin{equation}
\label{P1:asymp}
P_1(t) \simeq \frac{1}{\ln t}\,.
\end{equation}
\end{subequations}
This behavior strongly contrasts with $P_1(t)\simeq t^{-2}$ in critical
branching.  Because there is no mechanism for an isolated particle to die,
this leads to a much slower temporal decay of $P_1(t)$.

We now use this solution for $P_1(t)$ in conjunction with the rate equations
\eqref{master} and \eqref{single} to recursively find the asymptotic behavior
of $P_m(t)$.  Substituting \eqref{P1:asymp} into \eqref{single} leads to
\begin{equation*}
%\label{P2:asymp}
P_2(t) \simeq \frac{1}{2}\,\frac{1}{\ln t}
\end{equation*}
to leading order. Substituting these asymptotic forms for $P_1(t)$ and
$P_2(t)$ into the rate equation \eqref{master} with $m=2$, we find
\begin{equation*}
%\label{P3:asymp}
P_3(t) \simeq \frac{1}{3}\,\frac{1}{\ln t}
\end{equation*}
to leading order.  Continuing this procedure iteratively gives
\begin{equation}
\label{Pm:asymp}
P_m(t) \simeq \frac{1}{m}\,\frac{1}{\ln t}\,,
\end{equation}
which should be valid for $t\gg 1$ and $m\ll t$.

To determine the distribution $P_m(t)$ when $m\sim t$, it is helpful to first
determine the asymptotic growth of the moments $\langle m^k \rangle$.  
Using Eqs.~\eqref{master} and \eqref{single}, we deduce the rate equation $\frac{d \langle m \rangle}{d t} = P_1$. 
Combining this with the expression for $P_1(t)$ in \eqref{P1:asymp}, the
asymptotic growth of the average number of cells is
\begin{equation}
\label{m-av:growth}
\langle m \rangle \simeq \frac{t}{\ln t}\,.
\end{equation}
Similarly $\frac{d \langle m^2 \rangle}{d t} = 2 \langle m \rangle + P_1$ which in conjunction with \eqref{P1:asymp} and \eqref{m-av:growth} gives
\begin{equation}
\label{m2:growth}
\langle m^2 \rangle \simeq \frac{t^2}{\ln t}\,.
\end{equation}
For general $k\geq 2$, the evolution equation for $\langle m^k \rangle$ is
\begin{subequations}
\begin{equation}
\label{mk:growth-eq}
\frac{d \langle m^k \rangle}{d t} = 2 \sum_{a=1}^{\left\lfloor \frac{k}{2}\right\rfloor}\binom{k}{2a}\langle m^{k-2a+1} \rangle + P_1\,,
\end{equation}
whose leading behavior is
\begin{equation}
\label{mk:growth-asymp}
\frac{d \langle m^k \rangle}{d t} \simeq 2\binom{k}{2}\langle m^{k-1} \rangle\,.
\end{equation}
\end{subequations}
The results for the first two moments suggest that
$\langle m^k \rangle\simeq C_k\,t^k/\ln t$.  Substituting this ansatz into
\eqref{mk:growth-asymp}, we find that this ansatz is consistent and moreover
it serves to fix the amplitude $C_k$.  The final result is
\begin{equation}
\label{mk:growth}
\langle m^k \rangle \simeq (k-1)!\;\frac{t^k}{\ln t}
\end{equation}
for $k\geq 1$. (The zeroth moment is
$\langle m^0 \rangle= \sum_{m\geq 1}P_m(t)=1$ because extinction is
impossible.)

The small-$m$ tail of $P_m(t)$ and the behavior of the moments suggest that
the full distribution $P_m(t)$ has the following scaling behavior
\begin{equation}
\label{Pm:scaling}
P_m(t) \simeq \frac{1}{m}\;\frac{1}{\ln t}\;F(\mu), \qquad \mu = \frac{m}{t}\,,
\end{equation}
with $F(0)=1$. Using this distribution, we now express the moments
$\langle m^k \rangle$ through the integrals of the scaled distribution:
\begin{equation*}
%\label{mk:scaling}
\langle m^k \rangle \simeq \frac{t^k}{\ln t}\,\int_0^\infty d\mu\,\mu^{k-1} F(\mu)\,.
\end{equation*}
Thus we have the condition $\int_0^\infty d\mu\,\mu^{k-1} F(\mu) = (k-1)!$,
which allows us to fix the scaling function: $F(\mu) = e^{-\mu}$.  Finally,
we obtain
\begin{equation}
\label{Pm:scaling-final}
P_m(t) \simeq \frac{1}{m}\,\frac{1}{\ln t}\,e^{-m/t}\,.
\end{equation}
It is instructive to compare this scaling form with that of the critical
branching process.  From \eqref{classical}, this latter scaling form is
\begin{align}
  P_m(t)\simeq \frac{1}{t^2}\, e^{-m/t}\,.
\end{align}
We again see that the population distribution has a much slower temporal
decay than in critical branching; this occurs because of the effective ``source''
at $m=1$.

\section{Branching With Input}
\label{sec:stem}

A biologically motivated version of immortal branching is based on stem cells
that are always active.  We focus on a system that begins with a single stem
cell S.  This stem cell can give birth to ordinary (mortal) cells,
$\text{S}\to \text{S}+\text{P}$, and these mortal cells subsequently undergo
critical branching.  We denote the rate at which stem cells give birth as
$\beta$; for ordinary cells, the birth and death rates are set to 1.

Denote again by $P_m(t)$ the probability to have $m$ cells, i.e., the stem
cell and $m-1$ mortal cells.  The probability that only a stem cell exists obeys
\begin{subequations}
\label{stem:ME}
\begin{equation}
\label{single:stem}
\frac{dP_1}{dt} =  -\beta P_1+P_2\,,
\end{equation}
while for $m\geq 2$ the rate equations are
\begin{align}
\label{master:stem}
\frac{dP_m}{dt} &= m P_{m+1} -[2(m-1)+\beta]P_m +(m-2+\beta)P_{m-1}\,.
\end{align}
\end{subequations}

One can verify by direct substitution that
\begin{equation}
\label{Pm:stem}
P_m(t) =   \frac{\Gamma(m-1+\beta)}{\Gamma(m)\,\Gamma(\beta)}\;\frac{t^{m-1}}{(1+t)^{m-1+\beta}}
\end{equation}
is the solution of Eqs.~\eqref{stem:ME}. In the scaling region $m,t\to\infty$
with $\mu={m}/{t}$ finite, this solution acquires the scaling form
\begin{equation}
\label{Pmt:stem}
P_m(t) \simeq t^{-1}\Phi(\mu), \qquad \Phi(\mu)=\frac{\mu^{\beta-1}}{\Gamma(\beta)}\,e^{-\mu}\,.
\end{equation}
With this scaling form, the leading behavior of the moments are given by
\begin{equation}
\label{m-av-k:stem}
\langle m^k\rangle \simeq t^k\,\frac{\Gamma(\beta+k)}{\Gamma(\beta)}\,.
\end{equation}

The IBP described by Eqs.~\eqref{stem:ME} is often known as a branching
process with immigration.  The solution \eqref{Pm:stem} was derived by
Kendall~\cite{Kendall48}; see~\cite{Tibor} for recent work.  We now discuss a
multistage IBP that may have relevance to cancer, where tumor cells often
undergo multiple stages~\cite{armitage54,beer07}.  The simplest step in this
direction to allow for the possibility of two-type
branching~\cite{Kendall60}. These models are analytically
tractable~\cite{2-type}; models with three-type branching have not been
solved so far.

Two-type branching processes also give a good description of experimental
data~\cite{clayton07} on the cell dynamics of skin tissue.  The experimental
data are well fit by a model that involves progenitor cells (P) that can
divide (proliferate) and differentiate, and post-mitotic cells (M) that
eventually disappear (leave the basal layer).  The elemental steps of this
process are
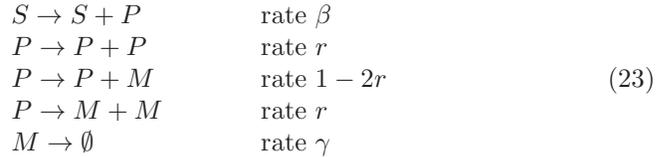
\begin{equation}
\label{full_def}
 \begin{tabular}{ll}
 $S \to  S+P$ \qquad\qquad & \quad rate $\beta$\\
 $P \to  P+P$                & \quad rate $r$\\
 $P \to  P+M$                &\quad rate $1-2r$\\
 $P \to  M+M$                &\quad rate $r$\\
 $M \to  \emptyset$       & \quad rate $\gamma$\\
 \end{tabular}
\end{equation}

With no source ($\beta=0$), this two-type branching process was solved in
\cite{skin}.  We now generalize the approach of Ref.~\cite{skin} to
$\beta>0$.  We account for the populations of progenitor and post-mitotic
cells by $P_{m,n}(t)$, the probability to have $m$ clones of $P$ and $n$
clones of $M$. This probability distribution evolves according to the master
equation
\begin{align}
\label{master_long}
 \frac{d P_{m,n}}{dt} &= [r(m-1)+\beta]P_{m-1,n} + r(m+1)P_{m+1,n-2} \nonumber\\
 & +  (1-2r)mP_{m,n-1} + \gamma (n+1)P_{m,n+1}  - (m+\gamma n + \beta)P_{m,n}\,.
\end{align}
We again introduce the  generating function 
\begin{equation*}
%\label{gener}
 \mathcal{P}(x,y,t) = \sum_{m,n=0}^\infty x^m y^n P_{m,n}(t)\,,
\end{equation*}
which satisfies the partial differential equation (PDE)
\begin{equation}
\label{P-eq}
\partial_t \mathcal{P} + U\partial_x \mathcal{P}+V \partial_y \mathcal{P} = W \mathcal{P} 
\end{equation}
where $U= x(1-y)- r(x-y)^2, ~V= \gamma (y-1)$ and $W= \beta(x-1)$.  For an
initially empty system $\mathcal{P}(x,y,t=0) = 1$.  Equation \eqref{P-eq} is
a first-order hyperbolic PDE and like all such equations, can be solved using
the method of characteristics~\cite{Logan08}. The characteristics are
determined by equations
\begin{align}
\label{char}
\begin{split}
   \frac{dx}{dt} &= U = x(1-y) - r(x-y)^2\\
   \frac{dy}{dt} &= V = \gamma(y-1)
 \end{split}
\end{align}
By construction of the characteristic equations, the time derivative 
along a characteristic is given by
\begin{equation*}
  \frac{d \mathcal{P}}{dt} = \partial_t \mathcal{P} + U\partial_x \mathcal{P}+V \partial_y \mathcal{P}
\end{equation*}
Combining this with Eq.~\eqref{P-eq} we get
$\frac{d \mathcal{P}}{dt} = W \mathcal{P}$, which we integrate to yield
\begin{equation}
\label{P-sol}
 \mathcal{P}(X,Y,T) = \exp\!\Big\{\beta\int_0^T d\tau\,[x(t)-1]\Big\}\,.
\end{equation}
It is convenient to think of $t$ as a `running' time, so we use the notation
$T$ for the final time.  The functions $x(t)$ and $y(t)$ vary along the
characteristic and rather than parameterizing by initial values $x(0)$ and
$y(0)$ it is preferable to parameterize by the final values $X=x(T)$ and
$Y=y(T)$.

To compute the generating function, Eq.~\eqref{P-sol}, we need to determine
$x(t)$ along the characteristic.  We integrate the second of Eqs.~\eqref{char}
to give
\begin{equation}
\label{y-eta}
y-1=(Y-1)e^{\gamma (t-T)}\,.
\end{equation}
Substituting this into the first of \eqref{char} and changing the time
variable $t$ to
\begin{equation}
\label{u-sol}
u=\gamma(1-Y)\,e^{\gamma (t-T)}\,,
\end{equation}
we arrive at 
\begin{equation}
\label{xu:eq}
 \frac{dx}{du} = \frac{x}{\gamma} - \frac{r}{\gamma u} (x-1+u)^2 \,.
\end{equation}

This Riccati equation admits an exact solution in terms of confluent
hypergeometric functions \cite{skin}.  Thus one can obtain a solution, but it
is complicated and difficult to use to obtain explicit results in the general
case.  A solution in terms of elementary functions is possible in the special
case $r=\frac{1}{4}$ and $\gamma=1$ (\ref{app:Riccati}).  Two prominent
results from this analysis are: (i) $P_m(t)$, the probability to have $m$
progenitor cells, and any number of post-mitotic cells and (ii) $P_{m,0}(t)$,
the probability to have $m$ progenitor cells, and \emph{no} post-mitotic
cells.  The former is given by
\begin{subequations}
\begin{equation}
\label{Pm-special}
P_m(t) = \frac{\Gamma(m+4\beta)}{\Gamma(4\beta)\,\Gamma(m+1)}\,
\frac{1}{(1+t/4)^{4\beta}}\, \left(\frac{t/4}{1+t/4}\right)^m \,.
\end{equation}
This result coincides with the exact solution \eqref{Pm:stem} with the
understanding that $m$ counts the number of progenitor cells, and the birth
and death rates are $r=1/4$.  The latter probability distribution is
\begin{equation}
\label{Pm0-special}
P_{m,0}(t) = \frac{\Gamma(m+4\beta)}{\Gamma(4\beta)\,\Gamma(m+1)}\,
\frac{e^{\beta(1-e^{-t})}}{(1+t/2)^{4\beta}}\, \left(\frac{t/4}{1+t/2}\right)^m \,.
\end{equation}
\end{subequations}

\section{Discussion}

We analyzed two immortal branching processes (IBPs). The first is a simple
extension of critical branching in which extinction cannot occur.  The
emergent behaviors are remarkably subtle, that are replete with logarithmic
corrections.  The second is a branching process with a steady input of cells,
or equivalently, immigration.  Such a steady input arises in a variety of
non-equilibrium many-body processes, such as aggregation with a steady input
of monomers~\cite{book}, fragmentation with a steady input of large
clusters~\cite{BK00}, and turbulence with a steady energy input at large
scales~\cite{B96,GS17}.  In these examples, many new features were uncovered
by incorporating steady input.  Rich behaviors also occur in spatially
extended systems with a spatially localized input~\cite{SR89,PK15} and the
role of spatial degrees of freedom in IBPs poses interesting challenges.

\section*{Acknowledgments}

We benefited from correspondence with Tibor Antal.
PLK thanks the hospitality of the Santa Fe Institute, where this work began.
The research of SR is partially supported by NSF grant DMR-1910736.

\appendix

\section{Derivation of Eq.~\eqref{Pm:final}}
\label{app:a}

To solve the governing differential equation \eqref{gen:eq} for the
generating function, we introduce the auxiliary variable $\zeta = 1/(1-z)$ to
recast this equation as
\begin{equation}
\label{gen:zeta}
s\mathcal{P} = \frac{d \mathcal{P}}{d\zeta} + 1 -\frac{\widetilde{P}_1(s) +1}{\zeta}\,.
\end{equation}
This linear inhomogeneous ordinary differential equation may be readily
solved to yield
\begin{equation}
\label{gen:sol}
\mathcal{P} = \frac{1}{s} - \big(\widetilde{P}_1(s) +1\big)
\int_0^\infty \frac{d\eta}{\eta+\zeta}\,e^{-s\eta}+ \Phi(s)\,e^{s\zeta}\,.
\end{equation}
To fix the integration constant $\Phi(s)$ we take the Laplace transform of
the normalization condition $\sum_{m\geq 1} P_m(t) = 1$ to give
$\sum_{m\geq 1} \widetilde{P}_m(s) = 1/s$.  This then implies that
$\lim_{z\to 1^-} \mathcal{P}(s,z) = 1/s$.  Since the limit $z\to 1^-$
corresponds to $\zeta\to +\infty$, the second term on the  right-hand side of
\eqref{gen:sol} vanishes while the third term, $\Phi(s)\,e^{s\zeta}$,
diverges if $\Phi(s)\ne 0$.  We thus conclude that $\Phi(s)=0$.

Returning to the original variable $z$, we rewrite \eqref{gen:sol} as
\begin{equation}
\label{gen:final}
\mathcal{P} = \frac{1}{s} - \big[\widetilde{P}_1(s) +1\big]\int_0^\infty d\eta\;e^{-s\eta}\, \frac{1-z}{1+\eta -\eta z}\,.
\end{equation}
We now expand in a Taylor series in $z$ to give the population distribution
\begin{equation}
\label{Pm:sol}
\widetilde{P}_m(s) = \big[\widetilde{P}_1(s) +1\big]\int_0^\infty d\eta\;e^{-s\eta}\,\frac{\eta^{m-1}}{(1+\eta)^{m+1}}\,.
\end{equation}
Specializing \eqref{Pm:sol} to $m=1$ we obtain a closed equation for
$\widetilde{P}_1(s)$. After some straightforward steps, this equation
simplifies to
\begin{equation}
\label{P1:sol}
\widetilde{P}_1(s) = \frac{1}{s e^s \Gamma(0,s)} - 1\,,
\end{equation}
where $\Gamma(0,s) = \int_s^\infty du\,e^{-u}/u$ is the incomplete
Gamma function.  Substituting the above form for $\widetilde{P}_1(s)$ in
\eqref{Pm:sol}, we obtain Eq.~\eqref{Pm:final}.

\section{Solution for $P_m(t)$ and $P_{m,0}(t)$ when $r=\frac{1}{4}$ and
  $\gamma=1$}
\label{app:Riccati}

For $r=\frac{1}{4}$ and $\gamma=1$, the Riccati equation \eqref{xu:eq} admits
a solution \cite{skin} in terms of elementary functions
\begin{equation}
\label{xu-sol}
x = 1+u+\frac{1}{\frac{1}{4}\ln u + C}\,.
\end{equation}
When $t=T$, we have $x=X$ and $u=1-Y$ for $\gamma=1$.  These requirements fix
the constant $C$ to be
\begin{equation}
C = - \tfrac{1}{4}\,\ln(1-Y) - (2-X-Y)^{-1}\,.
\end{equation}
Using this result we rewrite \eqref{xu-sol} as 
\begin{equation}
\label{xt-sol}
x(t) = 1+(1-Y)\,e^{t-T} - \left(\frac{T-t}{4} + \frac{1}{2-X-Y}\right)^{-1}\,.
\end{equation}
Substituting \eqref{xt-sol} into \eqref{P-sol}, computing the integral, and
returning to the original variables $x,y,t$, we obtain
\begin{equation}
\label{P-special}
 \mathcal{P}(x,y,t) = 
 \frac{\exp\!\big[\beta(1-e^{-t})(1-y)\big]}{\left[1+\frac{t}{4}(2-x-y)\right]^{4\beta}}\,.
\end{equation}

We now specialize \eqref{P-special} to the case $y=1$ and expand this
generating function in a Taylor series in $x$.  This gives
Eq.~\eqref{Pm-special}.  In a complementary way, we specialize
\eqref{P-special} to the case $x=1$ and use the definition
\begin{equation}
\label{gener-P}
 \mathcal{P}(1,y,t) = \sum_{n=0}^\infty y^n \Pi_{n}(t)\,,
\end{equation}
which expresses this restricted generating function in terms of $\Pi_{n}(t)$,
the probability to have $n$ post-mitotic cells.  Expanding
$\mathcal{P}(1,y,t)$ in a Taylor series in $y$ leads to
\begin{equation}
\label{Pn-special}
\Pi_n=e^{\beta(1-e^{-t})}\sum_{k=0}^n \frac{[-\beta(1-e^{-t})]^k}{k!}\,P_{n-k}(t) \,,
\end{equation}
with $P_m$ given by \eqref{Pm-special}. 

Using \eqref{P-special} we obtain a simple expression for the probability to
have no cells of any kind:
\begin{equation}
\label{P00}
P_{0,0}(t) =\mathcal{P}(0,0,t) = \frac{e^{\beta(1-e^{-t})}}{(1+t/2)^{4\beta}}\,.
\end{equation}
The temporal behavior of \eqref{P00} is the same as that for the probability
to have no progenitor cells, $P_0=(1+t/4)^{-4\beta}$
(Eq.~\eqref{Pm-special}), although with different amplitudes:
$P_0>P_{0,0}(t)$.  Finally, we specialize \eqref{P-special} to $y=0$ to give
\begin{equation}
\label{P-special-y=0}
 \mathcal{P}(x,0,t) = \sum_{m=0}^\infty x^m P_{m,0}(t)=
 \frac{\exp\!\big[\beta(1-e^{-t})\big]}{\left[1+\frac{t}{4}(2-x)\right]^{4\beta}}\,.
\end{equation}
Expanding this expression in a Taylor series gives the probabilities of
having only progenitor cells quoted in Eq.~\eqref{Pm0-special}.

\end{document}